\documentstyle[leqno,amssymb]{article}
 \def\qed{\unskip\quad \hbox{\vrule\vbox to 6pt {\hrule width
          4pt\vfill\hrule}\vrule} }
          \newcommand{\bez}{\nopagebreak\hspace*{\fill}
          \nolinebreak$\qed$\vspace{5mm}\par}

          \newenvironment{proof}{\vspace{1ex}
               \vspace*{10mm}\vspace{-10mm}\par\noindent{\bf
               Proof.}\nopagebreak
                    \par}{\nopagebreak\linebreak[0]\hspace*{\fill}
 $\qed$\vspace{5mm}\pagebreak[0]\vspace{3ex}\par}

                         \newtheorem{Th}{Theorem}

                         \newtheorem{Prop}{Proposition}

                         \newtheorem{Lemma}{Lemma}

                         \newtheorem{Remark}{Remark}

                         \newtheorem{Def}{Definition}

                         \newtheorem{Cor}{Corollary}

\def\a{\alpha}

\newcommand{\cb}{{\cal B}}

\newcommand{\ot}{\otimes}

\newcommand{\la}{\lambda}

\newcommand{\tf}{T_{\varphi}}

\newcommand{\va}{\varphi}

\newcommand{\vep}{\varepsilon}

\newcommand{\ov}{\overline}

\newcommand{\fqn}{\varphi^{(q_n)}}

\newcommand{\fpqn}{\left(\varphi^{(q_n)}\right)'}

\newcommand{\fbqn}{\overline{\varphi}^{(q_n)}}

\newcommand{\ftbqn}{\widetilde{\overline{\varphi}}^{(q_n)}}

\newcommand{\fpbqn}{\left(\overline{\varphi}^{(q_n)}\right)'}

\newcommand{\ftpbqn}{\left(\widetilde{\overline{\varphi}}^{(q_n)}\right)'}

\newcommand{\oinl}{\overline{I}_{n,l}}

\newcommand{\inl}{I_{n,l}}

\newcommand{\anl}{J_{n,l}(a,\varepsilon)}

\newcommand{\jnl}{J_{n,l}(a,\varepsilon)}

\newcommand{\beq}{\begin{equation}}

\newcommand{\eeq}{\end{equation}}

\newcommand{\Q}{{\Bbb{Q}}}
\newcommand{\R}{{\Bbb{R}}}

\newcommand{\T}{{\Bbb{T}}}

\newcommand{\Z}{{\Bbb{Z}}}

\newcommand{\N}{{\Bbb{N}}}

\newcommand{\xbm}{(X,{\cal B},\mu)}

\begin{document}

\title{On the ergodicity of cylindrical transformations given by the logarithm}
\author{Bassam Fayad \and
Mariusz Lema\'nczyk\footnote{Research partially supported by KBN
grant 1 P03A 03826.  \newline
}}

\maketitle \thispagestyle{empty}

\begin{abstract} 

Given $\a \in [0,1]$ and $\varphi: \T \to \R$ measurable,  the {\it cylindircal cascade} 
$S_{\a,\varphi}$ is the map from  $\T \times \R$ to itself  given by $S_{\a,\varphi} (x,y) = (x+\a,y+\varphi(x))$ that naturally appears in the study of some ordinary differential equations on $\R^3$. In this paper, we prove that for a set of full Lebesgue measure of $\a \in [0,1]$ the cylindrical cascades  $S_{\a,\varphi}$  are ergodic for
 every smooth function $\varphi$ with a logarithmic singularity, provided
 that the average of $\varphi$ vanishes.

Closely related to $S_{\a,\varphi}$ are the special flows  constructed above $R_\a$ and under $\varphi+c$ where  $c \in \R$ is such that $\varphi+c>0$. In the case of a function $\varphi$ with an  asymmetric logarithmic singularity our result gives the first examples of ergodic cascades $S_{\a,\varphi}$ with the corresponding special flows being mixing.
 Indeed, when the latter flows are mixing the usual techniques used to prove  the {\it essential value criterion} for $S_{\a,\varphi}$, that is equivalent to ergodicity, fail and we device a new method to prove this criterion that we hope could be useful in tackling other problems of ergodicity for cocycles preserving an infinite measure. 
  
\end{abstract}

\section{From flows to skew products}  Let  $(M,x_t,\nu)$
be a smooth dynamical system with continuous time and assume it
has a global section  $(\Sigma,T,\mu)$. For $\psi \in C^1(M,\R)$
one can consider the flow on $M \times \R$ given by coupling $x_t$
and the differential equation on $\R$
\begin{eqnarray} \label{flow}
{dz \over dt} = \psi(x_t), \quad z \in \R.
\end{eqnarray}
 The flow determined by the coupling has a skew product form and it
 is given by the formula
\beq\label{flow100} (x_0,z_0)\mapsto
(x_t,\int_0^t\psi(x_s)\,ds+z_0).\eeq It has also a section,
$\Sigma \times \R$, on which the dynamics writes as a skew product
over $T$, namely
\begin{eqnarray} \label{skew}
(\theta, z) \rightarrow (T\theta,z+\varphi(\theta)),
\end{eqnarray}
where $\varphi$ is obtained by integrating $\psi$ along flow
segments of $x_t$: $\varphi(\theta)=\int_0^t\psi(x_s)\,ds$, where
$t=t(\theta)$ is the first return time of $x_0=\theta$ to
$\Sigma$. In view of~(\ref{flow100})  the flow in (\ref{flow})
preserves the measure $\nu \times \lambda$, where $\lambda$
denotes Lebesgue measure on the line. When $(x_t,\nu)$, or
equivalently $(T,\mu)$, is ergodic, it is natural to ask whether
the flow given by (\ref{flow}) is ergodic for $\nu \times
\lambda$\footnote{Ergodicity for an infinite measure means that an
invariant set either has zero measure or its complement has zero
measure.}. This is equivalent to ergodicity of the skew product in
(\ref{skew}) for the measure $\mu \times \lambda$.

\begin{Remark} \em A necessary condition for ergodicity
of~(\ref{skew}) is that
 $\int_\Sigma \varphi(s)\, d \mu(s) =0$, which  by the Kac
 theorem we may always assume to hold
 by adding the constant
 $C=-\int_\Sigma\phi(\theta)\,d\mu/\int_{\Sigma}t(\theta)\,d\mu$ to $\psi$.
\end{Remark}

The study of skew products goes back to Poincar\'e and his work on
differential equations on $\R^3$ (see \S \ref{poincare} below
where $T$ is a minimal circular rotation and $\varphi$ is smooth)
and  was later undertaken in the general context, where on the
first coordinate, $T$ is an arbitrary ergodic automorphism of a
standard probability space $\xbm$, and on the second, $\varphi$ is
merely measurable (see monographs \cite{Aa} and \cite{Sch}).

In this note, we will prove the ergodicity of (\ref{skew}) when
$T$ is a minimal circular rotation $R_\a$, $\a$ belongs to a set
of full Lebesgue measure,  and $\varphi$ is a smooth function over
the circle except for an asymmetric logarithmic singularity (cf.\
the precise Definition~\ref{loglike} below). But first, we will
discuss the problems arising in the  study of the ergodicity of
(\ref{flow}) in the simplest case where $x_t$ is a smooth area
preserving flow on a surface and see how our result fits in this
context.

Note that when $x_t$ has only isolated fixed points od saddle
type, the global section $\Sigma$ exists and the return map $T$
will not be defined at the last points where $\Sigma$ intersects
the incoming separatrices of the fixed points and moreover the
return time function is asymptotic to infinity at these points.
Further, if $\psi$ does not vanish at a given fixed point, the
function $\varphi$ in (\ref{skew}) will have a singularity above
the corresponding point where $T$ is not defined and this
singularity will have the same nature as the one for the return
time function. It is not hard to see that a non-degenerate fixed
point of the saddle type of the flow $x_t$ yields a singularity of
the logarithmic type for the return time function.
 
\begin{Def} \label{loglike} \em We will say that a real function $\varphi$
defined over $\T$ has a {\em logarithmic singularity} at a point
$x_0$ if $\varphi$ is of class $C^2$ in $\T\setminus\{x_0\}$ and
there exist $A, B \in \R\setminus\{0\}$ such that
\begin{eqnarray*} \lim_{x \rightarrow x_0^-}  \varphi''(x) {(x-x_0)}^2 &=& A, \\
 \lim_{x \rightarrow x_0^+}  \varphi''(x) {(x-x_0)}^2 &=& B. \end{eqnarray*}
We say that the singularity is {\em asymmetric} if $A+B \neq 0$.
\end{Def}

\vspace{0.2cm}

\subsection{The case of linear flows on the torus}  \label{poincare}
When $x_t$ is an irrational flow on  the torus $\T^2$, it has a
global section $\T$ on which the Poincar\'e return map is a
minimal translation $R_\alpha$. The resulting skew products
\begin{equation} \label{skew2} S_{\a,\varphi}(\theta, z) = (\theta +\alpha,z+\varphi(\theta)),
\end{equation}
were intensively studied (for both $z \in \T$ and $z \in \R$) since they have been first introduced by Poincar\'e in \cite{Po}. 

Unlike the case $z \in \T$ where $S_{\a,\varphi}$ is ergodic (for the Haar measure of $\T^2$) if $\varphi$ equals a constant $\beta$ as soon as $1,\a,$ and $\beta$ are independent over $\Q$, a necessary condition for 
ergodicity in the case $z \in \R$ is that $\int_{\T} \varphi(\theta) d\theta =0$. In this case, the existence of ergodic skew products was first discovered by Krygin in \cite{Kr}. There exist elegant categorical proofs \cite{herman,katok} of the fact that the set of $(\a,\varphi)$ such that $S_{\a,\varphi}$ is ergodic forms a residual set (for the product topology) in the product of the circle with the space  $C^r_0(\T,\R)$ of functions of class $C^r$ with zero mean value (and this is true for any finite regularity $r \in \N$ or for $r=\infty$ or for the space $C^\omega_{\delta,0}(\T,\R)$  of real analtyic functions with zero mean value, analytically extendable in a fixed annular neighborhood of $\T$ of size $\delta$, continuous on its boundary, which is a Baire space if considered with the topology of uniform convergence).  Further, it actually holds that for a given Liouvillean $\a$, i.e. an $\a \in \R \setminus \Q$ such that 
$$\limsup_{p/q \in \Q} {-\log |\a - {p \over q}|  \over \log q}= \infty,$$
the set of $\varphi \in C^\infty_0(\T,\R)$ such that $S_{\a, \varphi}$ is ergodic is residual (for the $C^\infty$ topology),  and that for $\a$ satisfying 
$$\limsup_{p/q \in \Q} {-\log |\a-{p \over q}| \over q} \geq \delta>0,$$
then the  set of $\varphi \in C^\omega_{{\delta \over 2\pi},0}(\T,\R)$ such that $S_{\a, \varphi}$ is ergodic is residual (for the topology described above) (cf. e.g. \cite{Ba-Me}).

In specific situations however, proving ergodicity for skew products preserving an infinite measure may become a delicate task (cf. for example the problem of ergodicity raised in \cite{forrest}).  
 Ergodicity of $S_{\a,\varphi}$ was proved in  several situations,
e.g.: \cite{Aa-Le-Ma-Na}, \cite{Ba-Me}, \cite{Co}, \cite{Fr},
\cite{Kr}, \cite{Or}, \cite{Pa}, \cite{Vo}. \vspace{0.2cm}

\subsection{The case of time changed linear flows
on the torus with a stopping  point.} The easiest case of a flow
with a section where the Poincar\'e map is not defined at an
isolated point is a reparametrized irrational flow (multiply the
constant vector field by a smooth scalar function) on the torus
$\T^2$ where the orbit is stopped at an isolated point (isolated
zero for the reparametrizing function). But this procedure is not
interesting from the ergodic point of view because the flow thus
obtained is uniquely ergodic with respect to the Dirac measure
supported by the fixed point. The dynamics at the stopping point
is too slow (note that the inverse of the reparametrizing function
is not integrable, hence the flow preserves an infinite measure
which is equivalent to Lebesgue measure). This problem can be
bypassed by plugging in the phase space of the minimal linear flow
a weaker isolated singularity coming from a Hamiltonian flow in
$\R^2$. The so called Kochergin flows thus obtained preserve
beside the Dirac measure  at the singularity a measure that is
equivalent to Lebesgue measure. These flows still have $\T$ as a
global section with a minimal rotation for the return map, but the
slowing down near the fixed point produces a singularity for the
return time function above the last point where the section
intersects the incoming separatrix of the fixed point.  Again, if
$\psi$ does not vanish at the fixed point, this results in a
singularity of the same nature  for the function we obtain in the
system~(\ref{skew2}). The strength of the singularity depends on
how abruptly the linear flow is slowed down in the neighborhood of
the fixed point. A mild slowing down is typically represented by
the logarithm (e.g. \ when $\varphi(x)=-\log x-\log(1-x)-2$). In
this case ergodicity of (\ref{skew}) was proved in~\cite{Fr-Le}.
In the case of power like singularities, that were actually the
ones considered by Kochergin, no $\a \in \R \setminus \Q$ is known
for which we have ergodicity in~(\ref{skew}).

The second case is indeed sensitively different from the first one
for the following reason that we will further comment in the next
subsection: the special flow over $R_\a$ and under a smooth
function with at least one  power like singularity is
mixing~\cite{Koflows,bsmf} while the one under a smooth function
with symmetric logarithmic singularities is not~\cite{Ko,Le}.

\subsection{The case of a multi-valuated Hamiltonian on $\T^2$}
In \cite{Arnold}, Arnol'd investigated  Hamiltonian flows corresponding to
multi-valued Hamiltonians on a two dimensional torus for which the phase
space decomposes into cells that are filled up by periodic orbits and
one open ergodic component. On this component, the flow can be represented
as a special flow over a minimal rotation of the circle and under a ceiling
function that is smooth except for some logarithmic singularities.
The singularities are asymmetric since the coefficient in front of the
logarithm is twice as big on one side of the singularity as the one on the other side,
due to the existence of homoclinic saddle connections.

It follows that if $x_t$ in (\ref{flow}) is such a flow, the
system we obtain in~(\ref{skew}), once we restrict our attention
to the open ergodic component of $x_t$, is a skew product over a
minimal rotation of the circle with in the second coordinate a
function having asymmetric logarithmic singularities. In this
paper we prove the following.

\begin{Th}
\label{mmmain} For a.e.  $\alpha\in\T$, the cylindrical transformation $S_{\a,\varphi}: \T \times \R \to \T \times \R$, $(x,y)\mapsto(x+\alpha,\varphi(x)+y)$ is ergodic for any function
$\varphi$ of class $C^2$ on $\T \setminus \lbrace x_0 \rbrace$
with a logarithmic singularity at $x_0$ and with zero average.  \end{Th}

We do not know whether ergodicity holds for every irrational
$\alpha$, except for the special case when the singularity is
symmetric \cite{Fr-Le}.
Note that, unlike the symmetric case, the
special flows over irrational rotations and under smooth functions
with asymmetric logarithmic singularities are mixing
\cite{Kh-Si,Ko1,Ko2}. We will explain now why this fact makes the
usual proof of ergodicity of the skew product (\ref{skew2}) fail. 
We first need to introduce the {\it essential value criterion} which is necessary and sufficient  for the ergodicity of skew products.

Assume that $T$ is an ergodic automorphism of a standard
probability Borel space $\xbm$. Let $\va:X\to\R$ be a measurable
map. Denote by $\va^{(\cdot)}(\cdot):\Z\times X\to\R$ the cocycle
generated by $\va$, i.e. given by the formula
 \beq\label{e2} \va^{(n)}(x)=\left\{\begin{array}{ccc}
\va(x)+\va(Tx)+\ldots+\va(T^{n-1}x) &\mbox{if} & n>0\\
0 &\mbox{if}& n=0\\
-(\va(T^{n}x)+\ldots+\va(T^{-1}x))&\mbox{if}&
n<0\end{array}\right. \eeq  

Denote by
$\tf$ the transformation of $(X\times\R,\cb\ot\cb(\R),\mu\ot\la)$
given by
$$
\tf(x,y)=(Tx,\va(x)+y).$$ Note that
$(\tf)^n(x,y)=(T^nx,\va^{(n)}(x)+y)$ for each $n\in\Z$.

Following \cite{Sch} a number $a\in\R$ is called an {\em essential
value} of $\va$ if for each $A\in\cb$ of positive measure, for
each $\vep>0$ there exists $N\in\Z$ such that
$$
\mu(A\cap T^{-N}A\cap[|\va^{(N)}(\cdot)-a|<\vep])>0.$$ Denote by
$E(\va)$ the set of essential values of $\va$. Then the essential value criterion states as follows

\begin{Prop}\label{p1}
[\cite{Sch},\cite{Aa}]  We have

\noindent 1. $E(\va)$ is a closed subgroup of $\R$.

\noindent 2. $E(\va)=\R$ iff $\tf$ is ergodic.
\end{Prop}

Usual
methods of proving ergodicity of $S_{\a,\varphi}$ take into
consideration a sequence of distributions \beq\label{e3}
\left(\va^{(n_k)}\right)_\ast(\mu),\;k\geq1 \eeq (along some rigid
sequence $\{n_k\}$, i.e. $n_k\alpha\to0$ (mod 1) when
$k\to\infty$) as probability measures on the one-point
compactification of $\R$. As shown in \cite{Le-Pa-Vo} each point
in the topological support of a ``rigid" limit point of (\ref{e3})
is an essential value of the cocycle $\va$, hence contributing to
ergodicity of $S_{\a,\varphi}$. This method is especially well
adapted to those $\va$ whose Fourier transform satisfies
$\hat{\va}(n)=\mbox{O}(1/|n|)$, hence in particular for $\va$ of
bounded variation. The log symmetric $\va$ also enjoys this
property, see \cite{Fr-Le}, and indeed ergodicity in this case
holds over every irrational rotation. However the method fails in
the case of an asymmetric logarithmic function (or for functions
with power like singularities, no matter whether they are
symmetric or not)  since the distributions (\ref{e3}) tend to
Dirac measure at infinity.  The latter is indeed a necessary
condition for mixing of the corresponding special flows, cf.\
\cite{Le} or \cite{Sch1} for a more general case.
 
In the present note, in order to prove ergodicity of $\varphi$, we
will apply a different method which rather resembles
Aaaronson's abstract essential value condition (EVC) from
\cite{Aa-Le-Vo}.
 
To be more precise, the problem we face is the following: given $a \in \R$ and 
a rigidity sequence ${\{q_n\}}_{n \in \N}$ of $R_\a$, the sets $A_{n}(a,\epsilon)$ of points $x \in \T$  where
$\varphi^{(q_n)}(x) \in [a-\epsilon,a+\epsilon]$ have their measure tending to zero as $n$ goes to infinity; and if we ask that $q_n$ be a very strong rigidity sequence ($\a$ well approximated by rationals) so as to force $R_\a^{q_n}A_n(a,\epsilon)$ to self-intersect, we will not be able to have good lower bounds on the measure of the sets $A_n$ and it will be impossible therefore to show that $a$ is an essential value. If to the contrary we consider badly approximated numbers $\a$, $R_\a^{q_n} A_n(a,\epsilon)$ will be disjoint from $A_n(a,\epsilon)$ making the usual proof of the essential value fail. However, we stick to these numbers and  prove for some rigidity sequence
${\{q_n\}}_{n\in \N}$, that the sets $A_{n}(a,\epsilon)$ are not too
small (although their measure goes to zero),  i.e. that $\sum
\mu(A_{n}) = \infty$\footnote{This condition   fails  when we
consider functions with power like singularities and, in the case of asymmetric logarithmic
singularities, it holds only under some arithmetic restrictions 
 of Diophantine type on $\a$. For technical reasons, we do assume
however that,  along a sequence of integers with positive density,
the partial quotients of $\a$ are ``large enough".}, then we use
the structure of these sets on the circle and their almost
independence for different values of $n$ to deduce, using a
generalized version of the Borel-Cantelli lemma, that any measurable set can be measurably approximated by a union of $A_n$'s. We conclude after observing that the same holds for the sets $B_n = R_\a^{q_n}A_n$.

\subsection{Open problem: The general case of transitive  area preserving
flows with isolated singularities}  On surfaces of higher genus
the presence of fixed points is unavoidable for index reasons. For
area preserving flows with only isolated singularities, the return
map to any transversal is  conjugate to an interval exchange map.
Furthermore, if the flow is transitive then it is quasi-minimal,
i.e. every semi-orbit other than a fixed point or a point on a
separatrix of a saddle is dense. In general, the closure of any
transitive component is a surface with a quasi-minimal flow. If in
addition the fixed points are non-degenerate saddles then the
singularities of the return time function at the discontinuities
of the interval exchange map are of  logarithmic type. These
singularities are usually symmetric but asymmetric situations
similar to the one treated in the present paper  may appear, if
for instance there is a saddle point with one of its separatrices
forming a homoclinic saddle connection. In this general setting,
ergodicity of the underlying systems (\ref{flow}) is unknown:

\noindent {\bf Problem \ \ } Let $T:I \rightarrow I$ be an ergodic interval
exchange map. Let $\varphi$ be a smooth function defined over $I$ with
logarithmic singularities at the discontinuity points of $T$.
Assuming that $\int_I \varphi(\theta)d \theta =0$,
is $S: I \times \R \rightarrow I \times \R, (\theta,z)
\mapsto (T\theta,z + \varphi(\theta))$ ergodic?

\section{Notations. Properties of the sums \protect$\fqn$}

Throughout this text, $X$ will denote the additive circle $\T=\R/\Z$
 identified with [0,1) (mod 1). Recall (see e.g.\
\cite{Ch}) that each irrational number $\alpha\in[0,1)$ admits a
development into the continued fraction expansion
$$
\alpha=\frac{1}{\displaystyle{a_1+\frac{1}{a_2+\frac1{\ldots}}}},
$$ 
($a_i$ are positive integers) and $a_i$ are called the {\em
partial quotients} of $\alpha$, $i\geq 1$.  We have
$$
\frac{1}{2q_iq_{i+1}}<|\alpha-\frac{p_i}{q_i}|<\frac{1}{q_iq_{i+1}},
$$
where
$$
q_0=1, q_1=a_1, q_{i+1}=a_{i+1}q_i+q_{i-1}$$
$$
p_0=0,p_1=1, p_{i+1}=a_{i+1}p_i+p_{i-1}.$$ Recall also (e.g.
\cite{Ch}) that there exists a constant $c>1$ such that for $n$
large enough \beq\label{e4} q_n\geq c^n.\eeq

Untill the last section $\varphi$ will be
$$\va(x)=-1-\log(1-x), \quad x \in [0,1).$$
Note that $\va\in L^1(\T)$ and that $\int\va\,d\mu=0$.

If $f:\T\to\R$ is of bounded variations,  the following
Denjoy-Koksma inequality holds for the Birkhoff sums of $f$ along $R_\a$
$$
\left|\frac{1}{q_n}f^{(q_n)}(x)-\int_0^1
f\,d\mu\right| \leq \frac{1}{q_n}\mbox{Var}\,f
$$
for each $x\in[0,1)$ (see e.g. the proof of the Koksma inequality
in  \cite{Ku-Ni}).

Assume that $\alpha\in\T$ is irrational. Put
$$
H(\alpha)=\{n\geq0:\: q_{n+1}\geq 100
q_n\;\mbox{and}\;\alpha<\frac{p_n}{q_n}\}.
$$
Denote
$$
\oinl=\left[\frac{l}{q_n}+\frac{1}{50q_n},\frac{l+1}{q_n}-\frac{1}{50q_n}\right],
$$
$l=0,1,\ldots,q_n-1$.

\begin{Lemma}\label{l1} Assume that $H(\alpha)$ is infinite. Then for any
$a\in\R$, for all sufficiently large $n\in H(\alpha)$ we have:
\beq\label{f1} \mbox{$\fqn$ is continuous and strictly increasing
on each $\oinl$,} \eeq  \beq \label{f2} \left|\fpqn(x)-q_n\log
q_n\right|<\frac{1}{\sqrt{n}}q_n\log q_n\; \mbox{for
every}\;x\in\oinl, \eeq \beq\label{f3}
\fqn\left(\frac{l}{q_n}+\frac{3}{4q_n}\right)\geq a+1, \eeq
\beq\label{f4} \fqn\left(\frac{l}{q_n}+\frac{1}{4q_n}\right)\leq
a-1,\eeq $l=0,1,\ldots, q_n-1$.
\end{Lemma}
\begin{proof}
Denote
$$\ov{\va}(x)=\left(1-\chi_{[1-\frac{1}{50q_n},1]}(x)\right)\va(x),\;x\in[0,1).
$$
Assume that $n\in H(\alpha)$. We have
$$
\left|\alpha-\frac{p_n}{q_n}\right|\leq\frac{1}{100q_n^2}.
$$
Moreover, since $\alpha<\frac{p_n}{q_n}$, no point
$x,x+\alpha,\ldots,x+(q_n-1)\alpha$ belongs to
$[1-\frac{1}{50q_n},1)$ whenever $x\in\oinl$, $l=0,1,\ldots,
q_n-1$ (indeed,
$x+s\alpha=x+s\frac{p_n}{q_n}+s(\alpha-\frac{p_n}{q_n})$). It
follows that \beq\label{f1.1}
\fbqn(x)=\fqn(x)\;\;\mbox{for}\;\;x\in\bigcup_{l=0}^{q_n-1}\oinl.
\eeq Moreover, \beq\label{f1.2}
\mbox{Var}\,\ov{\va}=2\log(50q_n)-1. \eeq Integrating by parts the
integral $\int_0^{1-\frac{1}{50q_n}}\log(1-x)\,dx$ we find that
\beq\label{f1.3} \int_0^1\ov{\va}(x)\,dx=-\frac{\log q_n}{50q_n}.
\eeq We also have \beq\label{f1.4} \mbox{Var}\,\ov{\va}'=100q_n-1
\eeq and \beq\label{f1.5} \int_0^1\ov{\va}'(x)\,dx=\log(50q_n).
\eeq In view of (\ref{f1.1}) we have to show that the properties
(\ref{f1})-(\ref{f4}) hold for $\fbqn(x)$, $x\in\oinl$. Since no
point $x,x+\alpha,\ldots,x+(q_n-1)\alpha$
 belongs to $[1-\frac{1}{50q_n},1)$ and $\ov{\va}'$ is strictly positive
on $[0,1-\frac{1}{50q_n})$, (\ref{f1}) directly follows. Now, from
(\ref{f1.4}) and the Denjoy-Koksma inequality we obtain that
\beq\label{f1.7}
\left|\fpbqn(x)-q_n\int_0^1\ov{\va}'\,d\mu\right|\leq 100q_n-1.
\eeq Hence using (\ref{f1.5}) and (\ref{e4}), \beq\label{f1.8}
\left| \fpbqn(x)-q_n\log q_n\right|\leq\frac{1}{\sqrt{n}}q_n\log
q_n \eeq for $n$ large enough. Put
$$
\inl=\left[\frac{l}{q_n},\frac{l+1}{q_n}\right],
$$
$l=0,1,\ldots,q_n-1$ and
$$
\ftbqn(x)=\ov{\va}(x)+\ov{\va}(x+\frac{1}{q_n})+\ldots+\ov{\va}
(x+\frac{q_n-1}{q_n}),$$ $x\in[0,1)$. We have
$\int_{\inl}\ftbqn\,d\mu=\int_0^1\ov{\va}\,d\mu$, so by
(\ref{f1.5}), \beq\label{f1.9} \int_{\inl}\ftbqn\,d\mu=-\frac{\log
q_n}{50q_n}. \eeq In a similar manner as we proved (\ref{f1}) and
(\ref{f2}) we have that $\ftbqn$ is continuous and strictly
increasing on each $\inl$ and \beq\label{f1.10}
\left|\ftpbqn(x)-q_n\log q_n\right|<\frac{1}{\sqrt{n}}q_n\log q_n
\eeq for $n$ large enough. Moreover, \beq\label{f1.11}
\left|\fbqn(x)-\ftbqn(x)\right|\leq\frac{q_n\log
q_n}{q_{n+1}}\left(1+\frac{1}{\sqrt{n}}\right) \eeq for $n$ large
enough (and $x\in\oinl$). Indeed, for $x\in\oinl$, using the fact
that $\va'\geq0$ and that $i\frac{p_n}{q_n}>i\alpha$ for
$i=0,1,\ldots,q_n-1$, we have
$$
\left|\fbqn(x)-\ftbqn(x)\right|=\sum_{i=0}^{q_n-1}\ov{\va}'(\xi_{x,i})
\left(i\frac{p_n}{q_n}-i\alpha\right)$$ for some
$\xi_{x,i}\in\left[x+i\alpha,x+i\frac{p_n}{q_n}\right]$,
$i=0,1,\ldots, q_n-1$. Since
$0\leq\ov{\va}'(\xi_{x,i})\leq\ov{\va}'(x+i\frac{p_n}{q_n})$, we
obtain that
$$
\left|\fbqn(x)-\ftbqn(x)\right|\leq\frac{q_n}{q_nq_{n+1}}\sum_{i=0}^{q_n-1}
\ov{\va}'(x+i\frac{p_n}{q_n})=
$$
$$
\frac{1}{q_{n+1}}\ftpbqn(x)\leq\frac{q_n}{q_{n+1}}(1+\frac{1}{\sqrt{n}})\log
q_n
$$
and (\ref{f1.11}) follows.

In order to prove (\ref{f3}) it is hence enough to show that
\beq\label{f1.12}
\ftbqn\left(\frac{l}{q_n}+\frac{3}{4q_n}\right)\geq \frac{q_n\log
q_n}{q_{n+1}}\left(1+\frac{1}{\sqrt{n}}\right)+a+1. \eeq To show
(\ref{f1.12}), in view of (\ref{f1.10}) and the fact that
$q_{n+1}\geq 100q_n$, it is enough to show that
$$
\ftbqn\left(\frac{l}{q_n}+\left(\frac34-\frac15\right)\frac{1}{q_n}\right)\geq0
$$
(because the derivative of $\ftbqn$ is of order $q_n\log q_n$,
hence on the interval of length $\frac15\frac{1}{q_n}$ the
difference of the values of the function at the endpoints is at
least of order $q_n\log q_n\cdot\frac{1}{5q_n}=\frac15\log q_n$
which is bounded from below by the sequence of order
$\frac{q_n}{q_{n+1}}(1+\frac{1}{\sqrt{n}})\log q_n)$. Suppose to
the contrary that
$$
\ftbqn\left(\frac{l}{q_n}+\left(\frac34-\frac15\right)\frac{1}{q_n}\right)
\leq0.
$$
Using (\ref{f1.10}) consecutively for intervals
$\left[\frac{l}{q_n},\frac{l}{q_n}+(\frac34-\frac15)\frac{1}{q_n}\right]$
of length $(\frac34-\frac15)\frac{1}{q_n}$ and
$\left[\frac{l}{q_n}+\left(\frac34-\frac15\right)\frac{1}{q_n},
\frac{l+1}{q_n}\right]$ of length
$\left(\frac14+\frac15\right)\frac{1}{q_n}$ we find that
$$
\int_{\inl}\ftbqn\,d\mu\leq-\left(\frac34-\frac15\right)^2\frac{1}{q_n^2}
\left(1-\frac{1}{\sqrt{n}}\right)q_n\log q_n+
$$
$$
\left(\frac14+\frac15\right)^2\frac{1}{q_n^2}
\left(1+\frac{1}{\sqrt{n}}\right)q_n\log
q_n\leq-\frac1{11}\frac{\log q_n}{q_n},
$$
when $n$ is large enough, which is a contradiction with
(\ref{f1.9}).

In order to complete the proof it is enough to show that
$$
\ftbqn\left(\frac{l}{q_n}+\left(\frac14+\frac15\right)\frac{1}{q_n}\right)\leq
0.$$ Suppose the contrary. Then
$$
\int_{\inl}\ftbqn\,d\mu\geq\left(1-\frac{1}{\sqrt{n}}\right)
\left(\frac34-\frac15\right)^2 \frac{1}{q_n^2}q_n\log q_n-
$$
$$
\left(1+\frac{1}{\sqrt{n}}\right)\left(\frac14+\frac15\right)^2\frac{1}{q_n^2}q_n\log
q_n\geq 0$$ for $n$ large enough -- contradiction with
(\ref{f1.9}).
\end{proof}

\begin{Remark}\rm It is clear that small modifications in the proof
of Lemma~\ref{l1} will give us a similar result also in case
$\alpha>\frac{p_n}{q_n}$.
\end{Remark}

The lemma below will be essential in the proof of ergodicity of
$\va$.

\begin{Lemma}\label{l2}
For any $a\in\R$, any $0<\vep<1$, for any $n\in H(\alpha)$
sufficiently large there exists an interval
$$
\jnl\subset\left[\frac{l}{q_n}+\frac{1}{4q_n},\frac{l}{q_n}+\frac{3}
{4q_n}\right]$$ $(l=0,1,\ldots,q_n-1)$ such that for each
$x\in\jnl$, \beq\label{f2.1} \fqn(x)\in[a-\vep,a+\vep] \eeq and
\beq\label{f2.2} \left|\jnl\right|=\frac{2\vep}{q_n\log
q_n}+\mbox{o}\left(\frac{1}{q_n\log q_n}\right)\eeq.
\end{Lemma}
\begin{proof}
In view of (\ref{f1}), (\ref{f3}) and (\ref{f4}) of
Lemma~\ref{l1},
$$
\fqn\left(\left[\frac{l}{q_n}+\frac{1}{4q_n},
\frac{l}{q_n}+\frac{3}{4q_n}\right]\right)\subset[a-1,a+1],
$$
while the estimation (\ref{f2.2}) follows from (\ref{f2}).
\end{proof}

\section{Borel-Cantelli lemma and the Essential Value Criterion}
We will assume now that $\alpha$ satisfies: \beq\label{i1} n\in
H(\alpha)\;\;\mbox{for all}\;\; n\geq n_0, \eeq \beq\label{i2}
\sum_{n=1}^\infty\frac1{\log q_n}=+\infty. \eeq Fix $a\in\R$ and
$\vep>0$. Denote
$$
A_n=A_n(a,\vep)=\bigcup_{l=0}^{q_n-1}\anl.$$

\begin{Lemma}\label{l3}
For each $k\geq 1$,
$$
\sum_{n\geq k}\mu\left(A_n|\bigcap_{j=k}^{n-1}A_j^c\right)
=+\infty.$$
\end{Lemma}
\begin{proof}
First let us notice that the set $A_k^c$ is obtained from $[0,1)$
by discarding $q_k$ intervals $J_{k,l}(a,\vep)$, $l=0,1,\ldots,
q_k-1$, next the set $(A_k\cup J_{k+1})^c$ we obtain from $A_k^c$
by discarding $q_{k+1}$ intervals $J_{k+1,l}(a,\vep)$,
$l=0,1,\ldots, q_{k+1}-1$, and so on. At each step
$s=0,1,\ldots,n-1$ the set $\bigcap_{j=k}^{k+s}A_j^c$ is hence a
union of at most $q_k+q_{k+1}+\ldots+q_{k+s}+1$ consecutive,
pairwise disjoint intervals which we will call $s$-{\em holes}.
Call an $s$-hole {\em good} if its length is at least
$\frac{6}{q_{s+1}}$, otherwise it is called {\em bad}. Assume now
that $(a,b)$ is a good $s$-hole. At step $s+1$ we first divide
$[0,1)$ into $q_{s+1}$ intervals of equal length
$\frac1{q_{s+1}}$. Since $(a,b)$ is a good $s$-hole, we find
$$
0\leq r_1<r_2\leq q_{s+1}-1, r_2-r_1\geq 5\;\mbox{and}\;
\left[\frac{r_1+i}{q_{s+1}},\frac{r_1+i+1}{q_{s+1}}\right]\subset(a,b)$$
for each $i=0,1,\ldots, r_2-r_1-1$. We take $r_1$ and $r_2$
extremal with the above properties. For each $i=0,1,\ldots
r_2-r_1-1$ we then consider $J_{k+s+1,r_1+i}(a,\vep)$. We have
\beq\label{bassamprop}
J_{k+s+1,r_1+i}(a,\vep)\subset\left[\frac{r_1+i}{q_{s+1}}+\frac14\frac1{q_{s+1}},
\frac{r_1+i}{q_{s+1}}+\frac34\frac1{q_{s+1}}\right],\eeq
$i=0,1,\ldots, r_2-r_1-1$. Since $q_{n+1}\geq 100 q_n$, it follows
that $(a,b)$ is producing at least $r_2-r_1-1$ good $(s+1)$-holes.
Notice also that~(\ref{bassamprop}) and the inequality
$q_{n+1}\geq 100 q_n$ imply that any (either good or bad) $s$-hole
cannot produce more that two bad $(s+1)$-holes. With these
observations in hands we will show that for each $s\geq0$,
\beq\label{A} G_{k+s}\geq B_{k+s} \eeq where $G_{k+s}$ (resp.\
$B_{k+s}$) stands for the number of good (resp.\ bad) $s$-holes.
Indeed, for $s=0$, $B_{k+s}=0$. Assume that (\ref{A}) holds for
some $s\geq0$. Since each good $s$-hole produces at least
$r_2-r_1-1$ good $(s+1)$-holes, we have $G_{s+k+1}\geq 4 G_{k+s}$.
The number $B_{k+s+1}$ is bounded by $2G_{k+s}+2B_{k+s}$, whence
$G_{k+s+1}\geq B_{k+s+1}$ and (\ref{A}) follows.

Fix $s\geq0$ and consider the trace of $A_{k+s+1}$ on a good
$s$-hole $(a,b)$. There exists an absolute constant $c_1>0$ such
that
$$\mu\left(A_{k+s+1}\cap(a,b)\right)\geq c_1\mu(A_{k+s+1})\mu(a,b)
$$
(indeed, $\mu(A_{k+s+1})$ is of order $\frac{2\vep}{\log
q_{k+s+1}}$, $\mu(a,b)$ is of order $(r_2-r_1)\frac{1}{q_{k+s+1}}$
and $\mu(A_{k+s+1})\cap(a,b)$ is of order
$(r_2-r_1)\frac{2\vep}{q_{k+s+1}\log q_{k+s+1}}$). Taking into
account (\ref{A}), it follows that
$$
\mu\left(A_{k+s+1}\cap\bigcap_{j=0}^{s-1}A_{k+j}^c\right)\geq
\frac{c_1}2\mu(A_{k+s+1})\mu(\bigcap_{j=0}^{s-1}A_{k+j}^c).
$$
Hence
$$
\sum_{n\geq k}\mu(A_n|\bigcap_{j=k}^{n-1}A_j^c)\geq c_2\sum_{n\geq
k}\mu(A_n) \geq c_2\vep\sum_{n\geq k}\frac1{\log q_n}=+\infty
$$
and the lemma follows.\end{proof}  In what follows we will make use
of the following variant of the Borel-Cantelli lemma (see
\cite{Ne}, Prop. IV-4.4):

\vspace{3mm}

 Let $(\Omega,{\cal F},P)$ be a probability space. Let
$\{C_n\}\subset {\cal F}$. Suppose that for each $k\geq0$
$$
\sum_{n=k}^\infty P(C_n|\bigcap_{j=k}^{n-1} C_j^c)=+\infty.$$ Then
$$
\limsup_{n\to\infty}C_n=\Omega\;\;\mbox{(mod $P$)}.
$$

\vspace{3mm}

Directly from this and from Lemma \ref{l3} we obtain the
following.

\begin{Lemma}\label{l4} Under the above assumptions, $\limsup_{n\to\infty}
A_{n}(a,\vep)=\T$ (mod $\mu$).\bez
\end{Lemma}
Denote $B_n(a,\vep)=T^{q_n}A_n(a,\vep)$, $n\geq n_0$.

\begin{Lemma}\label{l5}
Under the above assumptions, $\limsup_{n\to\infty} B_n(a,\vep)=\T$
(mod $\mu$).
\end{Lemma}
\begin{proof}
Note that $T^{q_n}\anl$  is an interval of the same length as
$\anl$ and due to the condition
$|\alpha-\frac{p_n}{q_n}|<\frac1{100q_n^2}$ its position with
respect to $J_{n,l}(a,\vep)$ is not essentially changed. Therefore
we see that the arguments that lead to the proof of Lemma~\ref{l4}
work well also in case of the sequence $B_n(a,\vep)$, $n\geq n_0$.
\end{proof}

We are now able to prove that each real number is an essential
value of $\va$ under some restriction on $\alpha$.

\begin{Prop} \label{l6}
If $\alpha$ satisfies (\ref{i1}) and (\ref{i2}) then the logarithmic
cylindrical transformation is ergodic.
\end{Prop}
\begin{proof}
Take $a\in\R$. We will show that $a\in E(\va)$. Fix $0<\vep<1$. By
Lemmas~\ref{l4} and~\ref{l5}, for any $s\geq1$ we have (in
measure)
$$
\bigcup_{n=s}^\infty A_n=\T=\bigcup_{n=s}^\infty B_n,
$$
where $A_n=A_n(a,\vep)$, $B_n=B_n(a,\vep)$. Fix an interval $I$.
We have as $l$ goes to infinity, \beq\label{nul}
\mu(T^{q_l}I\triangle I)=\mu((I+q_l\alpha)\triangle I)\to 0. \eeq

Take an interval $\overline{I}$ that is strictly included in $I$
and such that $|\overline{I}| \geq \frac{99}{100} |I|$. For $s$
large enough the set $A_s = \bigcup_{n\geq s} \bigcup_{0\leq l
\leq q_n -1} J_{n,l} \cap \overline{I}$ satisfies $A_s \subset I$
and

\beq\label{g100} \mu(A_s)>\frac34|I|, \eeq likewise, using
(\ref{nul}), the set $B_s = \bigcup_{n\geq s} \bigcup_{0\leq l
\leq q_n -1} T^{q_n} J_{n,l} \cap \overline{I}$ satisfies $B_s
\subset I$ and

\beq\label{g101} \mu(B_s)>\frac34|I|. \eeq Note that if $x\in
A_s$, say $x\in \anl$, then $T^{q_n}x\in B_s\subset I$ and
$|\va^{(q_n)}(x)-a|<\vep$.

Finally, take  any Borel set $C\subset [0,1)$  of positive
measure. Let $x_0$ be a point of density. Take a small $\delta>0$
and let $I\ni x_0$ be an interval so that \beq\label{g102}
\mu(C\cap I)\geq (1-\delta)\mu(I). \eeq Taking into account
(\ref{g100}), (\ref{g101}) and (\ref{g102}), and choosing $\delta$
sufficiently small we obtain a pair $(n,l)$ such that
 the set
$$
\{x\in C:\: x\in J_{n,l} \;\mbox{and}\; T^{q_n}x\in  C\}
$$
is of positive measure and hence $a\in E(\va)$.
\end{proof}

\section{Proof of Theorem \ref{mmmain}}
In order to formulate the main result of this note, first notice
that to prove the assertion of Proposition~\ref{l6} we only need the
conditions (\ref{i1}) and (\ref{i2}) both to hold along a common
subsequence of denominators (indeed, in the proof of
Lemma~\ref{l3}, and hence of Lemmas~\ref{l4} and ~\ref{l5}, we
will consider the sets $\anl$ for $n$ belonging to the subsequence
and the relevant condition of independence needed to use the
Borel-Cantelli lemma also holds). Hence we have proved the
following.

\begin{Prop}\label{p2}
Assume that for $\alpha$ irrational there exists a subsequence
$\{n_k\}$ such that \beq\label{e24} q_{n_k+1}\geq 100 q_{n_k},\eeq
\beq\label{e25} \sum_{k=1}^\infty\frac1{\log q_{n_k}}=+\infty.\eeq
Then the cylindrical transformation
$(x,y)\mapsto(x+\alpha,-1-\log(1-x)+y)$ is ergodic. \bez\end{Prop}

Let us notice that the conditions (\ref{i2}) and~(\ref{e25}) are
almost equivalent in the following precise sense: (\ref{i2}) holds
if and only if (\ref{e25}) holds along an arbitrary subsequence
$\{n_k\}$ of positive lower density. (Indeed, positive lower
density of $\{n_k\}$ means that there exists a constant $M>0$ such
that $n_k\leq Mk$ for each $k\geq1$; write
$\{1,2,\ldots,Mn\}=\bigcup_k D_k$, where $D_k=\{kM,kM+1,\ldots,
(k+1)M-1\}$ and notice that given $k$, $\sum_{s\in D_k}\frac1{\log
q_{s}}\leq M\cdot\frac1{\log q_{n_k}}$ since the sequence
$\{\frac1{\log q_n}\}$ is decreasing.) The condition (\ref{i2}) is
satisfied for any $\alpha$ with bounded partial quotients. We
hence proved the following.

\begin{Cor}\label{c1}
Assume that $\alpha$ has bounded partial quotients. Assume that
there exists a subsequence $\{n_k\}$ of positive lower density
such that~(\ref{e24}) is satisfied along this subsequence. Then
the cylindrical logarithmic transformation is ergodic.\bez
\end{Cor}

\begin{Remark}\em
Let us notice that (inductively, using the formula
$q_{n+1}=a_{n+1}q_n+q_{n-1}$) we have
$$
a_1\ldots a_n\leq q_n\leq a_1\ldots a_n\cdot 2^n.$$ It follows
from this estimation that $$\sum_{n=1}^\infty \frac{1}{\log
q_n}=+\infty\;\mbox{iff}\; \sum_{n=1}^\infty
\frac{1}{\sum_{i=1}^n\log a_i}=+\infty.$$ Indeed, all we need  to
show is that
$$
\sum_{n=1}^\infty \frac{1}{\sum_{i=1}^n\log a_i}=+\infty
\;\mbox{iff}\; \sum_{n=1}^\infty \frac{1}{n+\sum_{i=1}^n\log
a_i}=+\infty.$$
This equivalence holds because as we have already noticed:\\
a series of positive decreasing frequencies is divergent iff it is
divergent
along a subsequence of positive lower density,
and moreover, given two increasing sequences a positive real
numbers, $\{b_n\}$, $\{c_n\}$ such that the series $\sum 1/b_n$
and $\sum 1/c_n$ diverge, also the series $\sum1/(b_n+c_n)$
diverges, for either on a set of positive lower density we have
$$
\frac{1}{b_n+c_n}\geq \frac{1}{2b_n} \;\mbox{or}\;
\frac{1}{b_n+c_n}\geq \frac{1}{2c_n}.$$
 \end{Remark}

\vspace{2mm}
 We claim now that the assumptions of
Proposition~\ref{p2} are satisfied for a.e.\ $\alpha\in\T$.
Indeed, we have that for a.e.\ irrational number $\alpha\in\T$,
$$
\lim_{n\to\infty}\frac{\log q_n}{n}=\frac{\pi^2}{12\log2}
$$
(see e.g. \cite{Co-Fo-Si}, Chapter 7), so the condition~(\ref{i2})
is satisfied for a.e.\ irrational $\alpha$. Then, consider the
Gauss transformation $x\mapsto Tx:=\{\frac1x\}$, $x\in(0,1)$ which
preserves the finite absolutely continuous measure
$dm=\frac1{1+x}\,dx$ with respect to which $T$ is mixing. We also
have $T^nx\in [1/(k+1),1/k)$ if and only if $a_n(x)=k$. Consider
$f(x)=\chi_{[1/(k+1),1/k)}(x)$. By the ergodic theorem, for a.e.\
$x\in(0,1)$,
$$
\lim_{N\to\infty}\frac1N\sum_{n=0}^{N-1}f(T^nx)=
\lim_{N\to\infty}\frac1N\sum_{n=0}^{N-1}f(a_n(x))=
m(\frac{1}{k+1},\frac1k)$$ and in particular the set of $n$'s such
that $a_n(x)=k$ has positive density. We hence proved
\begin{Prop}\label{p3}
The cylindrical transformation
$(x,y)\mapsto(x+\alpha,-1-\log(1-x) +y)$ is ergodic for a.e.
$\alpha\in\T$.\bez\end{Prop}

Note that all the calculations that were made for $\varphi(x) =
-1-\log(1-x)$ in view of Lemma \ref{l1} are  also valid for any
function of class $C^2$ on $\T \setminus \lbrace x_0 \rbrace$
having a logarithmic singularity at $x_0 \in \T$ (as in
Definition~\ref{loglike}) with $A=0$ and $B \neq 0$, and with zero
average.

Note also that Lemma \ref{l1} will hold for $\va_1=\va+f$ whenever
$f^{(q_n)}\to0$ uniformly, in particular when $f$ is absolutely
continuous and has zero mean (the uniform convergence to zero
follows from the Denjoy-Koksma inequality). Similarly, consider
the case of a function $\varphi_1$ having an asymmetric
logarithmic singularity at $0$. Then for some $D>0$ we  have that
$\varphi_1 = \tilde{\varphi} + f$ where $f(x)= - D \log x - D
\log(1-x), x \in (0,1)$, and $\tilde{\varphi}$ has a logarithmic
singularity at $0$ (as in definition~\ref{loglike}) with $A=0$ and
$B \neq 0$. Fix $0<\eta<1$ and let
$$
\ov{f}_n(x)=f(x)\cdot\chi_{\left[\eta/q_n,1-\eta/q_n\right]}.$$ We
have $\int_0^1\ov{f}_n'\,d\mu=0$, hence by the Denjoy-Koksma
inequality \beq\label{nons1} |(\ov{f}'_n)^{(q_n)}(x)|\leq
2q_n/\eta \eeq
 for each but finitely many $x\in\T$ (and for $n\geq n_0$). It follows
 that there exists a constant $c=c(\eta)$  such that if we put
 $$
 \tilde{I}_{n,l}=\left[\frac{l}{q_n}+ \frac{c}{q_n}, \frac{l+1}{q_n}-
 \frac{c}{q_n}\right]\;\;(l=0,1,\ldots,q_n-1)$$
then by the proof of Lemma~\ref{l1} we will
obtain~(\ref{f1})-(\ref{f4}) to hold on each $\tilde{I}_{n,l}$ if
we replace $\va$ by $\va_1$ and the RHS in the estimate (\ref{f2})
by $o(q_n\log q_n)$. It then follows that also Lemma~\ref{l2}
holds  and by repeating all the other other arguments we end up by
proving the following.

\begin{Th}\label{p4} For a.e.  $\alpha\in\T$,  the  cylindrical  transformation
$(x,y)\mapsto(x+\alpha,\varphi(x)+y)$ is ergodic for any function
$\varphi$ of class $C^2$ on $\T \setminus \lbrace x_0 \rbrace$
with an asymmetric logarithmic singularity at $x_0$ and with zero
average.  \bez \end{Th}

Theorem \ref{mmmain} then follows from this and the result of
\cite{Fr-Le} in the symmetric case.

\frenchspacing
\bibliographystyle{plain}

\vspace{5mm}

\noindent Laboratoire d'Analyse, G\'eom\'etrie
et Applications, UMR 7539\\
Universit\'e Paris 13 et CNRS\\
99, av. J.-B.\ Cl\'ement,\\
93430 Villetaneuse, France

\noindent fayadb@math.univ-paris13.fr

\vspace{5mm}

\noindent Faculty of Mathematics and Computer Science\\
Nicolaus Copernicus University\\
ul.\ Chopina 12/18\\
87-100 Toru\'n, Poland

\noindent mlem@mat.uni.torun.pl

\end{document}